# A Flight-Mechanics Solver for Aircraft Inverse Simulations and Application to 3D Mirage-III Maneuver

Osama A. Marzouk[*]

*University of Buraimi, College of Engineering, P.O. Box 890, P.C. 512, Al Buraimi, Sultanate of Oman*

**Abstract:** The main objective of this paper is to present a general mathematical model and an associated numerical algorithm applicable to an arbitrary fixed-wing fixed-mass aircraft undergoing an arbitrary maneuver, based on the 3D nonlinear coupled differential-algebraic equations of motion, including force, moment, kinematic and constraint equations. The model is formulated to address the inverse simulation problem where a target maneuver is prescribed and the corresponding time dependent patterns of the control variables are solved for to meet this maneuver. The model utilizes two different moving frames of references, namely the body axes and the wind axes. The numerical algorithm features sequential solution of equations in a fully explicit manner. It is straightforward to use the model in a reverse mode, namely the direct simulation problem.

The inverse problem may be summarized as follows:

**Inputs:** Time history of desired-trajectory rectangular coordinates relative to the ground-fixed axes. A constraint should be specified, which we arbitrarily chose it to be the bank angle. Also, certain geometric and aerodynamic aircraft data are needed.

**Outputs:** Time history of the control variables (thrust magnitude, elevator angle, rudder angle, ailerons angle), which will satisfy the aimed trajectory.

The paper finally applies the presented numerical algorithm to a roll maneuver for the Mirage-III fighter.

## 1. INTRODUCTION

The aircraft inverse simulation is a technique of great interest, where the equations of motion of an aircraft are solved for a given prescribed maneuver, leading to the prediction of needed time response of the control variables in order to achieve the required maneuver trajectory [1-5].

Great progress in the control and computer technology raised the interest in remotely piloted/operated vehicles (RPVs/ROVs), especially in the aerospace industry. Putting an unmanned control system on board of a flying vehicle helps in the scientific research, for example, without risking any human life. This can also be used for example to navigate a military RPV on a mission of photographing enemy sites in a certain area or to send an ultralight aircraft equipped with photographing cameras over a remote mission beyond visual observation range so direct radio control (RC) operation is not viable.

A building block for such programmed flight is an off-line flight mechanics simulator that can predict the necessary time-variation of the control variables to achieve a target trajectory. This scenario is called inverse simulation; and the reverse of which where the control variables are specified and the resultant trajectory is predicted is called direct simulation. The former is of practical compared to the latter.

The inverse simulation is also a useful design tool for aircrafts by allowing the prediction of the extreme values of the needed engine thrust force, deflection angles, and angle of attack for the expected maneuvers that the aircraft may perform. These values help the designer in sizing the engine (ensuring it can accommodate the needed thrust), determining the sufficient range of control surfaces deflections, and selecting a proper airfoil shape for the wing such that the experienced angle of attack does not reach the undesirable stall zone (which is accompanied with air separation from the upper wing surface and a sudden drop in the lifting force exerted on the wing). Being warned that the bare wing will face a stall problem, the designer can adopt one or more design changes to circumvent this, including the adoption of vortex generators [6].

We here present a mathematical model of flight mechanics for a generic fixed-wing fixed-mass aircraft undergoing a generic 3D (6 degrees of freedom, 6-DOF) motion described by a set of 18 nonlinear coupled differential-algebraic equations (DAE) in 18 variables. Out of the 18 variables, 4 are to be specified

[*]Address correspondence to this author at the University of Buraimi, College of Engineering, P.O. Box 890, P.C. 512, Al Buraimi, Sultanate of Oman; Tel: +968 91173908;
E-mail: osama.m@uob.edu.om

  

as known functions of time (or series of discrete values at certain time stations), and the rest of the variables are predicted by numerically solving the DAE system.

## 2. DERIVATION OF THE EQUATIONS OF MOTION

A 6-DOF simulation tracks the motion of a rigid body as it moves through the atmosphere. The moving vehicle can rotate as well as translate. The 6 fundamental differential equations that allow the motion to be tracked are presented below, derived from applying Newton's second law for linear and angular momentums. The body-fixed reference frame is utilized in the derivation.

### 2.1. Three Force Equations (Linear Momentum)

In formulating these equations, the weight force (mg) is handled in an explicit form separate from the X, Y, Z force components; which account for all other forces acting on the aircraft. The weight acts in the direction of the ground axis $z_g$.

If the translational equations are formulated in the body axes, we obtain:

$$m \left( \frac{dV}{dt}\bigg|_{body} + \omega \times [V]_{body} \right) = [F]_{body} \quad (a)$$

where the subscript (body) indicates that the time derivative or the vector components are with respect to the moving body axes. The above vector equation translates into the three following scalar equations [7]:

$$m(\dot{u} - vr + wq) = X$$
$$m(\dot{v} - wp + ur) = Y \quad (b)$$
$$m(\dot{w} - uq + vp) = Z$$

While Equations (a) and (b) can be employed in theory to solve the equation of motion, they may suffer from poor efficiency from the numerical point of view due to the potential large disparity in the size of the adjacent terms. To demonstrate this, consider the case of a supersonic aircraft has a flight speed of 600 m/s having a reasonable upper limit on pitch-rate q of about 2 rad/s. In this case, the artificial acceleration term (u q) can be as large as 1200 m/s$^2$ or 122 g's. On the other hand, the actual acceleration term (Z/m), the acceleration due to the external force in $z_b$, primarily weight and aerodynamic lift, may have an upper limit of only a few g's. Thus, the artificial accelerations can be greater than the actual accelerations by a factor of say 50 due to the high rotation rates, which the body-axes experience. This means unfavorable computer scaling and thus poorer solution accuracy for a given computer precision. In addition, Equations (a) and (b) couple the high-speed dynamics of the rotational equations into the translational equations; this places severe computational demands.

To resolve the above problems, the translational equations are formulated with reference to the wind (flight-path) axes. In doing this, the following geometrical relations are implicitly used (see Figures **3–5**):

$$\tan \alpha = w/u$$
$$\sin \beta = v/V \quad (c)$$
$$\cos \alpha \cos \beta = u/V$$

#### 2.1.1. $x_w$ Component

$$T \cos \alpha \cos \beta = -\bar{q} S \begin{pmatrix} C_x \cos \alpha \cos \beta + \\ C_y \sin \beta + C_z \sin \alpha \cos \beta \end{pmatrix}$$
$$- mg \, (\cos \theta \sin \phi \sin \beta - \sin \theta \cos \alpha \quad (1)$$
$$\cos \beta - \cos \theta \cos \phi \sin \alpha \cos \beta) + m\dot{V}$$

#### 2.1.2. $y_w$ Component

$$mV\dot{\beta} = mg \, (\cos \theta \sin \phi \cos \beta + \sin \theta \cos \alpha \sin \beta) - T \cos \alpha \sin \beta$$
$$+ \bar{q} S \, (C_y \cos \beta - C_x \cos \alpha \sin \beta - C_z \sin \alpha \sin \beta) \quad (2)$$
$$+ mV(-r \cos \alpha + p \sin \alpha) - mg \sin \alpha \cos \theta \cos \phi \sin \beta$$

#### 2.1.3. $z_w$ Component

$$mV\dot{\alpha} \cos \beta = mg \, (\sin \theta \sin \alpha + \cos \theta \cos \phi \cos \alpha)$$
$$+ \bar{q} S C_z \cos \alpha - (T + \bar{q} S C_x) \sin \alpha \quad (3)$$
$$+ mV(q \cos \beta - r \sin \beta \sin \alpha - p \sin \beta \cos \alpha)$$

### 2.2. Three Moment Equations (Angular Momentum)

The rotational equations of motion should be formulated in the body axes, as done below. This is because in this case the moments of inertia will be constants rather than being time-varying quantities.

#### 2.2.1. $x_b$ Component

$$T_0 \dot{p} = (BC - D^2) T_1 + (FC + ED) T_2 + (FD - EB) T_3 \quad (4)$$

#### 2.2.2. $y_b$ Component

$$T_0 \dot{q} = (FC + ED) T_1 + (AC - E^2) T_2 + (AD + EF) T_3 \quad (5)$$

#### 2.2.3. $z_b$ Component

$$T_0 \dot{r} = (FD + BE) T_1 + (AD + FE) T_2 + (AB - F^2) T_3 \quad (6)$$

where

$$T_0 = A\,B\,C - A\,D^2 - B\,E^2 - C\,F^2 - 2\,D\,E\,F$$
$$T_1 = (B-C)\,q\,r + (E\,q - F\,r)\,p + (q^2 - r^2)\,D + L$$
$$T_2 = (C-A)\,r\,p + (F\,r - D\,p)\,q + (r^2 - p^2)\,E + M$$
$$T_3 = (A-B)\,p\,q + (D\,p - E\,q)\,r + (p^2 - q^2)\,F + N$$

The above equations can be simplified if the aircraft exhibits a geometric symmetry about the $x_b$–$z_b$ plane, in which case D ($I_{yz}$) = 0 and F ($I_{xy}$) = 0.

## 2.3. Three Equations Relating the Euler Angle Rates and Body Angular Rates

These equations relate the Euler rates $(\dot{\phi}, \dot{\theta}, \dot{\psi})$ with the angular velocities referenced to the body axes (i.e., p, q, r).

As the simulated body translates and rotates due to the external forces and moments, its attitude must be traced in order to determine its true orientation.

The simplest approach computes the Euler angle rates from the three body-fixed rates. These latter rates are merely the results of numerically integrating the body rotational accelerations given by Equations (4,5,6). The mathematical relationship between the Euler angle rates and the body-fixed rotational rates is now presented.

### 2.3.1. Euler Angles

Let the angular velocity vector ω of a moving body be expressed in terms of its components in the body-fixed axes,

$$\omega = p\,\hat{e}_{xb} + q\,\hat{e}_{yb} + r\,\hat{e}_{zb} \quad (d)$$

The angular velocity can also be expressed in terms of its components along the three axes of rotation relating the local-level and the body-fixed coordinate systems, as shown in Figure **1**. The local coordinate system is a moving rectangular coordinate system with the aircraft, sharing the same origin with the body-fixed system; however, the local system is always aligned with the ground system.

$$\omega = \dot{\phi}\,\hat{e}_{xb} + \dot{\theta}\,\hat{e}_{y'} + \dot{\psi}\,\hat{e}_{zL} \quad (e)$$

The local-level unit vectors along $z_L$ and along intermediate y´ are expressed in terms of the body-fixed system unit vectors as:

$$\hat{e}_{zL} = -\sin\theta\,\hat{e}_{xb} + \cos\theta\sin\phi\,\hat{e}_{yb} + \cos\theta\cos\phi\,\hat{e}_{zb} \quad (f)$$

$$\hat{e}_{y'} = \cos\phi\,\hat{e}_{yb} - \sin\phi\,\hat{e}_{zb} \quad (g)$$

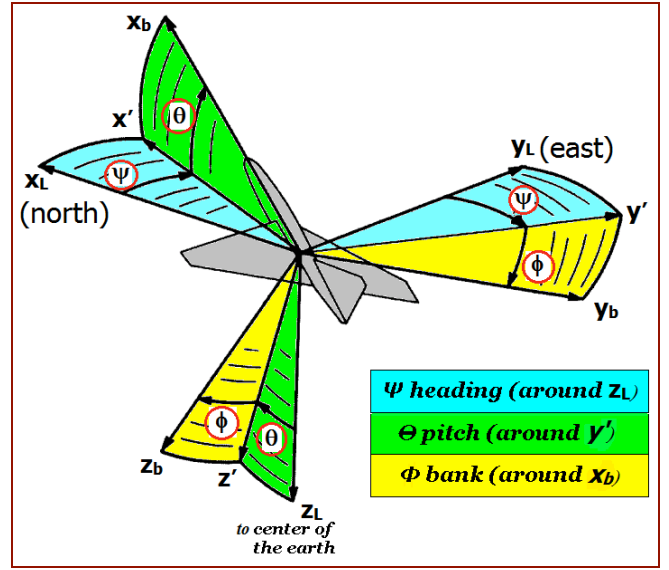

**Figure 1:** Relationship between local-level and body-fixed coordinate systems (adapted from [8]).

Equating Equations (d) and (e) for ω, and substituting with Equations (f) and (g), the Euler angular rates can be expressed in terms of the body-fixed rates as follows:

$$p = -\sin\theta\,\dot{\psi} + \dot{\phi} \quad (7)$$

$$q = \sin\phi\cos\theta\,\dot{\psi} + \cos\phi\,\dot{\theta} \quad (8)$$

$$r = \cos\phi\cos\theta\,\dot{\psi} - \sin\phi\,\dot{\theta} \quad (9)$$

## 2.4. Three Equations Relating the Velocity Components with the Flight Path Angles

The following expressions derive from straight-forward resolution of the effective velocity vector V into its components along the ground axes (Figure **8**).

$$\dot{x}_g = V\cos\theta_w\cos\psi_w \quad (10)$$

$$\dot{y}_g = V\cos\theta_w\sin\psi_w \quad (11)$$

$$\dot{z}_g = -V\sin\theta_w \quad (12)$$

## 2.5. Two Equations Relating Flight Attitude Angles with the Flight Path Angles

$$\cos\theta_w\sin(\psi_w - \psi) = \sin\beta\cos\phi - \cos\beta\sin\alpha\sin\phi \quad (13)$$

$$\sin\theta_w = \cos\beta\cos\alpha\sin\theta - (\sin\beta\sin\phi + \cos\beta\sin\alpha\cos\phi)\cos\theta \quad (14)$$

It is worthy of being mentioned that in the special case where the velocity vector is aligned with the longitudinal axis (thus, α = β = 0), then the flight attitude angles (ψ, θ) become identical to the flight path angles ($\psi_w$, $\theta_w$). In that case, Equation (13) becomes a trivial equality of zeros, and Equation (14) reduces to a equality between $\sin(\theta_w)$ and $\sin(\theta)$. However, that specialty does not stipulate relational condition among the flight attitude angles (ψ, θ, φ).

## 2.6. Aerodynamic Forces and Moments

There are 3 external forces and 3 external moments that act on the simulated rigid-body aircraft. These forces and moments are most easily handled in the body-fixed coordinate system.

The forces and moments are functions of nondimensional aerodynamic coefficients which should be available either based on numerical prediction using appropriate software packages (like the Advanced Airplane Analysis [9]) or based on wind tunnel tests.

The coefficients are converted to forces using the dynamic pressure and a reference area (being the projected area of the wing, as typically done in aeronautical applications [10]); and further converted to moments using an additional reference length.

The following expressions describe a basic set of aerodynamic forces and moments that act on the aircraft:

$X \equiv F_X = C_x\, \bar{q}\, S$　　aerodynamic force in body axis $x_b$

$Y \equiv F_Y = Cy\, \bar{q}\, S$　　aerodynamic force in body axis $y_b$

$Z \equiv F_Z = C_z\, \bar{q}\, S$　　aerodynamic force in body axis $z_b$

$L \equiv M_X = C_l\, \bar{q}\, S\, d$　　aerodynamic moment about body axis $x_b$

$M \equiv M_Y = C_m\, \bar{q}\, S\, d$　　aerodynamic moment about body axis $y_b$

$N \equiv M_Z = C_n\, \bar{q}\, S\, d$　　aerodynamic moment about body axis $z_b$

### 2.6.1. Aerodynamic and Stability Coefficients

The nondimensional lift coefficient is approximated as a linear function of the angle of attack [11]. It should be noted that while this simple equation is often used by aeronautical engineers, it breaks down remarkably at a relatively high α, around 15°, due to the separation of the boundary layer from the wing skin, leading to a sharp decline in the lift coefficient. For many airfoil sections of wings, the $C_L$ has a non-zero value (denoted by $C_{L0}$) at α = 0. This is due to the asymmetry (camber) of the airfoil section. For symmetric airfoils, $C_{L0}$ is 0. We consider the general asymmetric airfoil having

$$C_L = C_{L0} + C_{L\alpha}\, \alpha \quad (15)$$

where $C_{L\alpha} = \dfrac{dC_L}{d\alpha}$. The nondimensional drag coefficient is most often expressed as a quadratic function of the lift coefficient, following the so-called drag polar [12] which is

$$C_D = C_{D0} + K_{CD}\, C_L^2 \quad (16)$$

The two previous equations are rewritten as:

$$C_D = C_{D0} + \frac{C_{D\alpha}}{2C_{L\alpha}}\, C_L \; ; \; C_{D\alpha} = 2\, K_{CD} C_L C_{L\alpha} \quad (17)$$

We also introduce the side force coefficient, $C_C$

$$C_C = C_{Y?}\, \beta \quad (18)$$

With this, we have now

$$\begin{aligned}
C_x &= -C_D \cos\alpha \cos\beta - C_C \cos\alpha \sin\beta + C_L \sin\alpha \\
C_y &= -C_D \sin\beta + C_C \cos\beta \\
C_z &= -C_D \sin\alpha \cos\beta - C_C \sin\alpha \sin\beta - C_L \cos\alpha \\
C_l &= C_{l\beta}\, \beta + C_{lp}\, (p\, b/v) + C_{lr}\, (r\, b/v) + C_{l\delta l}\, \delta_l + C_{l\delta n}\, \delta_n \\
C_n &= C_{n\beta}\, \beta + C_{np}\, (p\, b/v) + C_{nr}\, (r\, b/v) + C_{n\delta l}\, \delta_l + C_{n\delta n}\, \delta_n \\
C_m &= C_{m0} + C_{m\alpha}\, \alpha + C_{mq}\, q + C_{m\delta m}\, \delta_m
\end{aligned} \quad (19)$$

**Notes:**

- From the above expressions, it is pointed out that the control-surface displacements appear implicitly in the equations of angular momentum, Equations (4,5,6).

- In the very special case of level steady flight without acceleration or rotation, the DAE system reduces to only two non-trivial equations, namely: $T = \tfrac{1}{2} \rho\, S\, V^2 C_D$ from Equation (1) and $0 = mg - \tfrac{1}{2} \rho\, S\, V^2 C_L$ from Equation (3). These These are the two fundamental equilibrium equations for cruising (steady and horizontal) flight, setting a force balance between the thrust and drag, and between the lift and weight.

## 3. AUGMENTED DAE SYSTEM

As a summary of the presented mathematical model, we have 14 differential-algebraic equations formed by Equations (1-14) that must be solved to get the 18 following variables:

$V$, $\alpha$, $\beta$

$p$, $q$, $r$

$\theta$, $\phi$, $\psi$

$\delta_l$, $\delta_m$, $\delta_n$, $T$

$\psi_W$, $\theta_W$

$\dot{x}_g$, $\dot{y}_g$, $\dot{z}_g$

Therefore, we must have an additional 18 -14 = 4 constraint equations, i.e. 4 of these variables should be known a priori in time, either as analytical functions of time or as sequences of sampled values.

### 3.1. Constraint Equations

In order to solve the equations of motion, 4 constraint equations are augmented into the differential-algebraic equations system, expanding it into 18 non-linear, coupled equations in 18 variables. The 4 constraint equations are

$x_g = fun_1(t)$

$y_g = fun_2(t)$

$z_g = fun_3(t)$

$\phi = fun_4(t)$ or one of the 4 controls will be arbitrarily chosen (e.g. T = constant).

In the procedure used here for solving the system of DAEs, we will take the 4$^{th}$ constraint to be $\phi = fun_4(t)$, i.e. the bank angle is known. So the main output from solving the system is the four control variables $\delta_l$, $\delta_m$, $\delta_n$, $T$.

## 4. SOLUTION PROCEDURE

It should be noted that the trajectory period will be discretized into a certain number of points (time stations), separated by a uniform time step (Δt) and the solution will be obtained at these stations. We need additional relational expressions to algebraically evaluate the time derivatives of some variables at each time station, and these expressions are obtained by differentiating some equations in the DAE system with respect to time. All the equations in the 18-equation DAE system (thus the equation to be differentiated here) are functions of some of the 18 variables appeared in that system. The solution procedure here solves the system sequentially as a series of scalar equations without the need of a solver for a system of equations.

### 4.1. Extra Derivative Equations

- From Equations (4,5,6), it is shown that values of $\dot{p}, \dot{q}, \dot{r}$ are needed in order to obtain the control displacements $\delta_l$, $\delta_m$, $\delta_n$. Expressions for $\dot{p}, \dot{q}, \dot{r}$ are to be obtained from differentiating Equations (7,8,9) once with respect to time. The resulting expressions for $\dot{p}, \dot{q}, \dot{r}$ will be referred as Equations ($\dot{7}, \dot{8}, \dot{9}$), respectively.

- Then $\ddot{\theta}, \ddot{\phi}, \ddot{\psi}$ will appear in the right-hand side of Equations ($\dot{7}, \dot{8}, \dot{9}$). Thus, expressions for them are needed. Because $\ddot{\phi}$ should be available from processing the 4$^{th}$ constraint, then we are left with $\ddot{\psi}$ and $\ddot{\theta}$. Twice differentiating Equations (13 and 14) yields expressions for $\ddot{\psi}$ and $\ddot{\theta}$, and we refer to them as Equations ($\ddot{13}$ and $\ddot{14}$), respectively. We will also need the first time derivates $\dot{\psi}$ and $\dot{\theta}$ to initialize the solution (see subsection 4.4. Initialization), and thus will keep the expressions for the first derivative of Equations (13 and 14) and denote them by Equations ($\dot{13}$ and $\dot{14}$), respectively.

- Now, $\ddot{\beta}$ and $\ddot{\alpha}$ appear in the right-hand side of Equations ($\ddot{13}$ and $\ddot{14}$). Needed expressions for them will be obtained from differentiating once Equations (2 and 3), which we refer to them as Equations ($\dot{2}$ and $\dot{3}$), respectively.

- Then, $\dot{T}$ appears in the right-hand side of Equations ($\dot{2}$ and $\dot{3}$). An expression for it is to be obtained from differentiating once Equation (1), yielding Equation ($\dot{1}$).

- Differentiating Equations (11,12) once and twice with respect to time, it becomes possible to evaluate algebraically the values of $\dot{\theta}, \dot{\psi}_W, \ddot{\theta}_W, \ddot{\psi}_W$ at all time stations if we know the corresponding values and derivates of the ground-based coordinates $y_g$ and $z_g$ and the velocity $V$. These resulting expressions are referred to as Equations ($\dot{11}, \dot{12}, \ddot{11}, \ddot{12}$).

### 4.2. Phases of the Solution

We divide the solution procedure into three phases:

### 4.3. Setup

1. The ground-based coordinates $x_g$, $y_g$, $z_g$ are given from the trajectory input functions (fun$_1$, fun$_2$, and

fun$_3$). Calculate their time derivatives up to the third order.

2. Use Equations (10,11,12) to solve algebraically for the values of $V, \theta_W, \psi_W$ at all time stations.

3. Use Equations $(1\dot{1}, 1\dot{2}, 1\ddot{1}, 1\ddot{2})$ to evaluate values of $\dot{\theta}, \dot{\psi}_W, \ddot{\theta}_W, \ddot{\psi}_W$ at all time stations.

4. Use finite difference expressions to evaluate $\dot{V}$ and $\ddot{V}$ at all stations from the obtained values of V.

For the terminal time stations, either forward or backward (one-sided) differences are used, but central double-sided difference is used for general stations. The following second-order-accurate expressions are used:

$$\dot{V}_n = \frac{1}{2\Delta t}(-3V_n + 4V_{n+1} - V_{n+2})$$
$$\ddot{V}_n = \frac{1}{\Delta t^2}(2V_n - 5V_{n+1} + 4V_{n+2} - V_{n+3})$$
$\quad \ldots \quad n = 1$ (initial time station) (20)

$$\dot{V}_n = \frac{1}{2\Delta t}(3V_n - 4V_{n-1} + V_{n-2})$$
$$\ddot{V}_n = \frac{1}{\Delta t^2}(2V_n - 5V_{n-1} + 4V_{n-2} - V_{n-3})$$
$\quad \ldots \quad n = nMax$ (final time station) (21)

$$\dot{V}_n = \frac{1}{2\Delta t}(V_{n+1} - V_{n-1})$$
$$\ddot{V}_n = \frac{1}{\Delta t^2}(V_{n+1} - 2V_n + V_{n-1})$$
$\quad \ldots \quad 1 < n < nMax$ (general time station) (22)

These expressions can be derived by imagining 3 consecutive values of V to which a quadratic function in the form V(t´) is fitted, where (t´) is a local time variable being zero at the nth station. Differentiating this fitting function once (for $\dot{V}$) or twice (for $\ddot{V}$) and evaluating the result at t´ = 0 gives the above finite difference expressions. For the expressions involving 4 points, a cubic fit is needed (a quadratic fit is insufficient to achieve second-order accuracy) [13].

5. The bank angle $\phi$ is given as input function (fun4). Calculate its first and second derivatives at each time station.

### 4.4. Initialization

The first (initial) time station is assigned the time station index n = 1.

1. The angles $\beta_{ini}$ and $\alpha_{ini}$ are initially equal to zeros.

2. The angles $\psi_{ini}$ and $\theta_{ini}$ are calculated from Equations (13 and 14).

3. The initial thrust magnitude $T_{ini}$ is calculated from Equation (1).

4. The initial time derivates $\dot{\beta}_{ini}$ and $\dot{\alpha}_{ini}$ are zeros.

5. The initial time derivates $\dot{\psi}_{ini}$ and $\dot{\theta}_{ini}$ are calculated from the expressions in Equations ($1\dot{3}$, $1\dot{4}$).

6. The initial body-fixed angular rates $p_{ini}, q_{ini}, r_{ini}$ are evaluated using Equations (7,8,9), respectively.

7. The initial time derivates $\dot{p}_{ini}, \dot{q}_{ini}, \dot{r}_{ini}$ are set to zeros.

8. The initial deflection angles $\delta_{l,ini}, \delta_{m,ini}, \delta_{n,ini}$ are calculated from Equations (4,5,6).

#### 4.4.1. About $C_{L0}$ and α

It is worth mentioning here that the above initial values satisfy equilibrium. So the angle of attack (α) in the procedure refers actually to the change in angle of attack from the equilibrium value. So when we obtain from the programmed procedure that α=0, this actually means that $\alpha_{actual} - \alpha_{equb} = 0$, where the subscript (equb) refers to the equilibrium condition: lift force ($\bar{q}$ S $C_L$) = weight (m g). A good treatment (used here) to streamline this issue within the calculations is to set the (0) reference condition in $C_{L0}$ to be at the equilibrium point rather than at the α=0 point. With this, we have

$$C_{L0,procedure} = C_{L0,equb} = \frac{mg}{\bar{q}S}$$
$$\alpha_{actual} = \alpha_{procedure} + \alpha_{equb} - |\alpha_{zero\text{-}lift}|$$
$$\alpha_{equb} = \frac{C_{L0,equb}}{C_{L\alpha}}$$
$$\alpha_{zero\text{-}lift} = -\frac{C_{L0,actual}}{C_{L\alpha}}$$
(23)

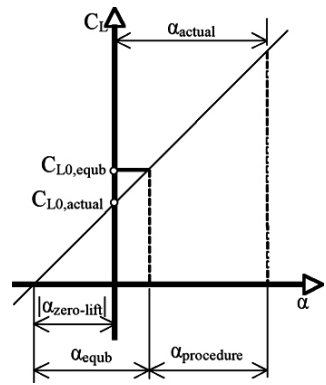

**Figure 2:** Relation between the angle of attack in the procedure and the actual one.

The relation among the different values of the he angle of attack is illustrated graphically in Figure **2**. In fact, this issue will not affect any of our results except the thrust.

### 4.4.2. Atmospheric Air Density

The following relationship is used to describe the dependence of the air density on the altitude (which is equal to $-z_g$). This expression is derived for the thermal-gradient troposphere layer of the atmosphere, which extends from the sea level to an altitude of about 11 km.

$$\rho = 1.225\left(1 + \frac{0.0065}{288} z_g\right)^{\left(1 - \frac{g}{0.0065 * 287}\right)} \quad (24)$$

The value of 1.225 is the standard air density at sea level in kg/m$^3$. The 0.0065 is the lapse rate (rate of temperature decrease with altitude, being about 6.5 °C/km [14]). The 288 is the standard air temperature at sea level in kelvins. The 287 in the exponent is the specific gas constant of air, in J/kg.K. In the above expression, the altitude is in meters and g = 9.81 m/s$^2$.

### 4.5. Loop of Calculation

The old values of the 12 variables $\left(\alpha,\beta,\theta,\psi,T,\dot{\alpha},\dot{\beta},\dot{\theta},\dot{\psi},p,q,r\right)_n$ for the old time station (n) should be known, either from the initial values (when n = 2) or from the results of the previous station, indexed: n-1 (when n > 2).

The core of this phase is a fourth-order Runge-Kutta (RK4) integration method, which can be summarized as:

For a generic ordinary differential equation of the form: $\dot{y} = f(t)$, subject to a starting condition of $y_{old}(t_{old})$, the RK4 method increments to a new value of $y_{new}(t_{old} + \Delta t)$ through

$$\left(k\dot{y}\right)_1 = f(y_0) \quad ; \quad y_{temp,1} = y_{old} + \left(k\dot{y}\right)_1 \Delta t / 2$$

$$\left(k\dot{y}\right)_2 = f(y_{temp,1}) \quad ; \quad y_{temp,2} = y_{old} + \left(k\dot{y}\right)_2 \Delta t / 2$$

$$\left(k\dot{y}\right)_3 = f(y_{temp,2}) \quad ; \quad y_{temp,3} = y_{old} + \left(k\dot{y}\right)_3 \Delta t$$

$$\left(k\dot{y}\right)_4 = f(y_{temp,3})$$

$$y_{new}(t_{old} + \Delta t) = y_{old} + \Delta t \left(\frac{\left(k\dot{y}\right)_1 + 2\left(k\dot{y}\right)_2 + 2\left(k\dot{y}\right)_3 + \left(k\dot{y}\right)_4}{6}\right)$$

Now, going back to our computational problem, the RK4 method is used as follows:

1. Fourth-order Runge-Kutta loop ($k_j$ ; j =1,2,3,4)

   i. $k\dot{T}$ is calculated from Equation(1).

   ii. $k\ddot{\beta}$ and $k\ddot{\alpha}$ are calculated Equations ($\dot{2}$ and $\dot{3}$), respectively.

   iii. 
   $$\left(k\dot{\beta}\right)_1 = \dot{\beta}_{old} \qquad\qquad \left(k\dot{\alpha}\right)_1 = \dot{\alpha}_{old}$$
   $$\left(k\dot{\beta}\right)_2 = \dot{\beta}_{old} + \left(k\ddot{\beta}\right)_1 \Delta t / 2 \quad \left(k\dot{\alpha}\right)_2 = \dot{\alpha}_{old} + \left(k\ddot{\alpha}\right)_1 \Delta t / 2$$
   $$\left(k\dot{\beta}\right)_3 = \dot{\beta}_{old} + \left(k\ddot{\beta}\right)_2 \Delta t / 2 \quad \left(k\dot{\alpha}\right)_3 = \dot{\alpha}_{old} + \left(k\ddot{\alpha}\right)_2 \Delta t / 2$$
   $$\left(k\dot{\beta}\right)_4 = \dot{\beta}_{old} + \left(k\ddot{\beta}\right)_3 \Delta t \quad \left(k\dot{\alpha}\right)_4 = \dot{\alpha}_{old} + \left(k\ddot{\alpha}\right)_3 \Delta t$$

   iv. $k\ddot{\psi}$ and $k\ddot{\theta}$ are calculated from Equations ($1\dot{3}$ and $1\dot{4}$), respectively.

   v. 
   $$\left(k\dot{\psi}\right)_1 = \dot{\psi}_{old} \qquad\qquad \left(k\dot{\theta}\right)_1 = \dot{\theta}_{old}$$
   $$\left(k\dot{\psi}\right)_2 = \dot{\psi}_{old} + \left(k\ddot{\psi}\right)_1 \Delta t / 2 \quad \left(k\dot{\theta}\right)_2 = \dot{\theta}_{old} + \left(k\ddot{\theta}\right)_1 \Delta t / 2$$
   $$\left(k\dot{\psi}\right)_3 = \dot{\psi}_{old} + \left(k\ddot{\psi}\right)_2 \Delta t / 2 \quad \left(k\dot{\theta}\right)_3 = \dot{\theta}_{old} + \left(k\ddot{\theta}\right)_2 \Delta t / 2$$
   $$\left(k\dot{\psi}\right)_4 = \dot{\psi}_{old} + \left(k\ddot{\psi}\right)_3 \Delta t \quad \left(k\dot{\theta}\right)_4 = \dot{\theta}_{old} + \left(k\ddot{\theta}\right)_3 \Delta t$$

   vi. $k\dot{p}, k\dot{q}, k\dot{r}$ are calculated from Equations ($\dot{7},\dot{8},\dot{9}$), respectively.

2. From the respective values of ($k_1$, $k_2$, $k_3$, $k_4$) for each of the 12 variables, evaluate $\left(\alpha,\beta,\theta,\psi,T,\dot{\alpha},\dot{\beta},\dot{\theta},\dot{\psi},p,q,r\right)_{n+1}$ at the new time station (n+1).

3. Use these new values at the new time station: (n+1) to calculate algebraically $\dot{T}_{n+1}$ from Equation ($\dot{1}$)

4. Calculate algebraically $\left(\ddot{\beta},\ddot{\alpha}\right)_{n+1}$ from Equations ($\dot{2},\dot{3}$)

5. Calculate algebraically $\left(\ddot{\psi},\ddot{\theta}\right)_{n+1}$ from Equations ($1\dot{3},1\dot{4}$)

6. Calculate algebraically $\left(\dot{p},\dot{q},\dot{r}\right)_{n+1}$ from Equations ($\dot{7},\dot{8},\dot{9}$) or as the weighted average of the four values of $k\dot{p}, k\dot{q}, k\dot{r}$. The latter choice is easier and very reasonable, especially for a small time interval ($\Delta t$), and it is implemented in our work – for example: $\dot{p}_{n+1} = \left[\left(k\dot{p}\right)_1 + 2\left(k\dot{p}\right)_2 + 2\left(k\dot{p}\right)_3 + \left(k\dot{p}\right)_4\right]/6$.

7. From $\left(\dot{p},\dot{q},\dot{r}\right)_{n+1}$, and with $\left(V,\alpha,\beta,p,q,r\right)_{n+1}$ known, the new control deflection angles $\left(\delta_l, \delta_m, \delta_n\right)_{n+1}$ are obtained from Equations (4,5,6).

8. Repeat steps (1–7) for each subsequent time station.

## 5. TERMINOLOGY

### 5.1. Ground Axes

Imaginary coordinate system ($x_g$, $y_g$, $z_g$) fixed on earth, by which we can describe variables with respect to the earth. The ground $x_g$–$y_g$ plane is parallel to the horizon. We can think of $x_g$ as arbitrarily pointing to the geographic north; accordingly, $y_g$ is pointing to the geographic east. Although the earth is rotating, this system is considered to be inertial (Newtonian, non-accelerating) and the effects of such rotation may be safely neglected for aeronautical applications because the induced acceleration due to this rotation is relatively small [15, 16].

### 5.2. Body Axes

Another coordinate system ($x_b$, $y_b$, $z_b$), fixed on the aircraft with its origin located at the center of gravity of the aircraft (as shown in Figure **3**). The body-fixed axis $x_b$ points toward the nose of the aircraft. It coincides with the fuselage centerline. The body-fixed axis $y_b$ points toward the starboard side of the aircraft (the right side of the pilot).

The body axis $x_b$ is also called the longitudinal axis or the roll axis. The body axis $y_b$ is also called the lateral axis or the pitch axis. The body axis $z_b$ axis is also called the yaw axis.

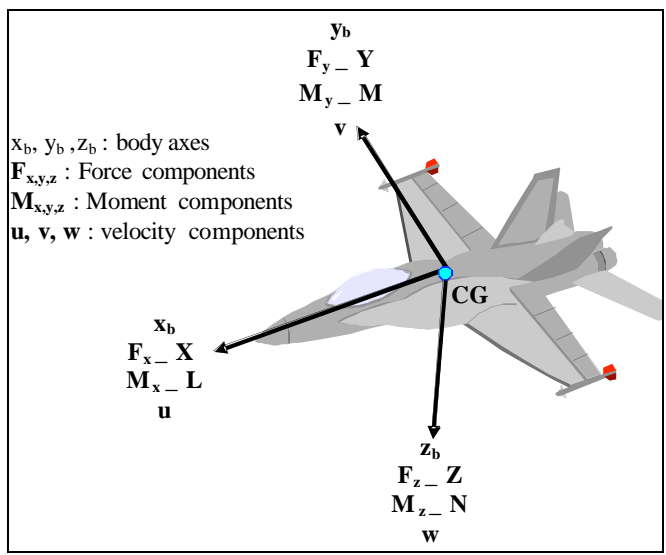

**Figure 3:** Body axes (and the relevant components of the force, moment, and effective/relative velocity).

### 5.3. Wind Axes

Another coordinate system ($x_w$, $y_w$, $z_w$) moving with aircraft. It considers the aircraft as a point particle located at its center of gravity and pays no attention to the aircraft orientation.

The wind axes system aligns its $x_w$ with the velocity vector of the aircraft's center of gravity relative to the wind. It is thus tangent to the aircraft path. In this paper, the term 'velocity' refers by default to this relative velocity because it is the effective velocity as far as the aerodynamic loads exerted on the aircraft are concerned. The direction of the $y_w$ is uniquely decided such that it is aligned with the ground axis $y_g$ when the $x_w$ is aligned with the ground axis $x_g$.

### 5.4. Trajectory

A path described by the time history of the coordinates measured from the origin of the ground axes, along which an aircraft can fly.

### 5.5. Maneuver

The ability of an aircraft to fly on a certain path using the various controls on the aircraft. These controls include the deflection angles of the moving surfaces (like rudder) and the thrust.

## 6. ENGLISH SYMBOLS (AND CORRESPONDING SI UNITS)

| | |
|---|---|
| ***A, B, C*** | Moments of inertia about rolling, pitching, and yawing axes respectively. They are also denoted elsewhere by *$I_{xx}$, $I_{yy}$, $I_{zz}$*, respectively; **kg.m$^2$** |
| ***b*** | Reference length for lateral derivatives (such as the wing span); **m** |
| ***d*** | Reference length for longitudinal derivatives (such as the wing mean chord); **m** |
| ***CG*** | Center of gravity (center of mass) of the aircraft |
| ***$C_L$, $C_D$, $C_M$*** | Lift, drag, and pitching-moment coefficients, respectively; **dimensionless** |
| ***D, E, F*** | Products of inertia in the body-axes planes $y_b$-$z_b$, $z_b$-$x_b$, and $x_b$-$y_b$, respectively. They are also denoted elsewhere by *$I_{yz}$, $I_{zx}$, $I_{xy}$*, respectively; **kg.m$^2$** |
| $\hat{e}$ | Unit vector; **dimensionless** |
| ***g*** | Gravitational acceleration; constant at 9.81 **m/s$^2$** |
| **L, M, N** | Rolling, pitching, and yawing moments, respectively. They are also denoted by $M_x$, $M_y$, $M_z$, respectively; **N.m** |
| **m** | Mass of the aircraft; **kg** |
| **p, q, r** | Angular velocities of the rolling, pitching, and yawing, respectively; **rad/s** |
| $\bar{q}$ | Dynamic pressure $= \frac{1}{2}\rho V^2$; **N/m²** |
| **S** | Wing projected area; **m²** |
| **V** | Magnitude of the wind/aircraft relative velocity; **m/s** |
| **T** | Magnitude of the thrust, which is directed toward the longitudinal axis; **N** |
| **t** | Time ; **s** |
| **u, v, w** | Components (longitudinal, lateral, and downward, respectively) of the relative velocity vector in the body axes; **m/s** |
| | As per particle curvilinear kinematics, the velocity vector is instantaneously tangent to the flight path made by the aircraft CG. |
| **X, Y, Z** | Components of the total applied force, except the weight, in the body axes. They are also denoted by $F_x$, $F_y$, $F_z$, respectively; **N** |

## 7. GREEK SYMBOLS

| | |
|---|---|
| $\alpha$ | Angle of attack (AOA), Figure **4**; rad |
| | The angle of attack is the angle that the velocity vector (*V*) makes with the body-axes horizontal plane ($x_b$-$y_b$), being positive when *V* goes toward $z_b$. |

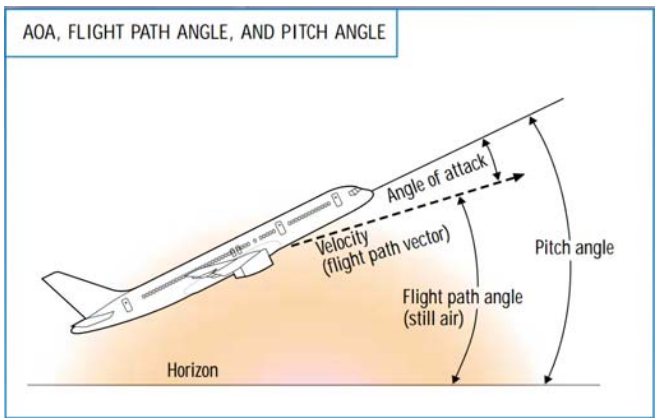

**Figure 4:** Explanation of the angle of attack (AOA: α), flight path angle ($\theta_w$), and pitch angle (θ) [17].

| | |
|---|---|
| $\beta$ | Sideslip angle, which is the angle that the velocity vector (*V*) makes with the body-axes plane of symmetry ($x_b$-$z_b$), being positive when *V* goes toward $y_b$ (Figure **5**); **rad** |

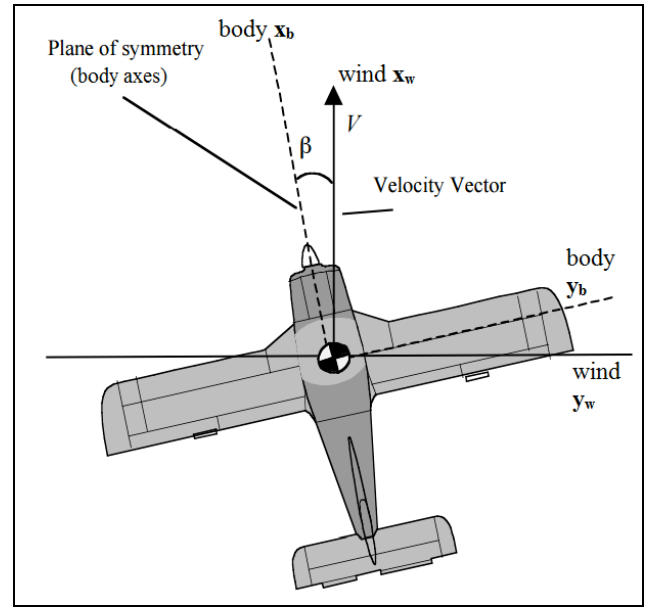

**Figure 5:** Illustration of the sideslip angle (adapted from [18]).

| | |
|---|---|
| $\delta_l, \delta_m, \delta_n$ | Deflections angles of the movable control surfaces: aileron, elevator, and rudder, respectively; **rad** |
| | The positive sense of these angles is depicted in Figure **6**. |
| | $\delta_l$: aileron deflection is positive if the trailing edge of right/starboard aileron moves up |
| | $\delta_m$: elevator deflection is positive if the trailing edge moves down |
| | $\delta_n$: rudder deflection is positive if the trailing edge moves toward the left of the pilot |

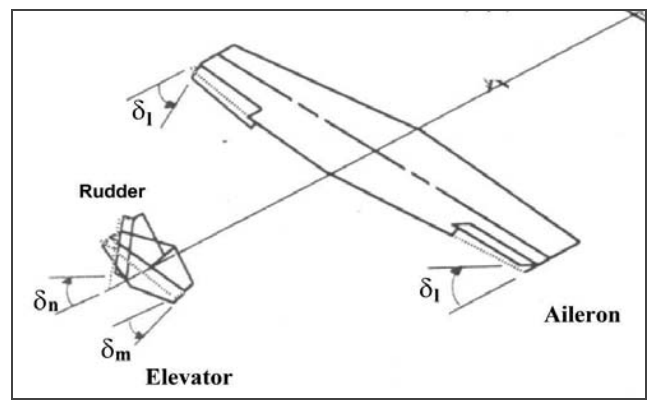

**Figure 6:** Positive directions of the deflection angles of the control surfaces.

$\psi, \theta, \phi$     Flight attitude angles: bank angle (also called roll angle), heading angle (also called yaw angle), and pitch angle, respectively (see Figure **7**); **rad**

These angle are also called Euler angles (see subsection Euler Angles), Tait–Bryan Euler angles, and Tait–Bryan angles. They have the following ranges:

$$-\pi \leq \psi \leq \pi$$
$$-\tfrac{\pi}{2} \leq \theta \leq \tfrac{\pi}{2}$$
$$-\pi \leq \phi \leq \pi$$

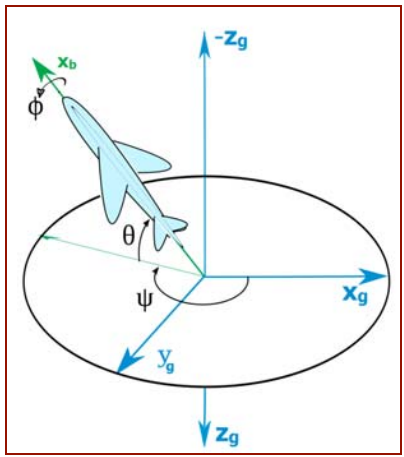

**Figure 7:** Illustration of the flight attitude angles (adapted from [19]).

$\psi_w, \theta_w$     Flight path angles. They describe the direction of the velocity vector *V* with respect to the ground axes (see Figure **8**); **rad**

The angle $\theta_w$ is the angle of *V* above the horizon. The angle $\psi_w$ is the angle of *V* relative to the global $x_g$-$z_g$ plane, being positive when *V* goes towards $y_g$.

$\rho$     Air density; **kg/m³**

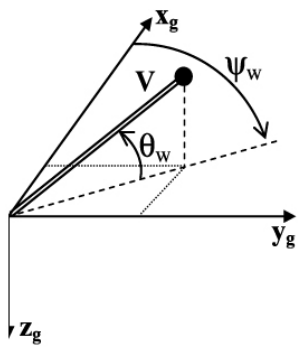

**Figure 8:** Illustration of the flight path angles.

## 8. RESULTS AND DISCUSSION

We implemented the computational procedure for computer execution as a MATLAB[20] m-script file and then applied it to a test case. This case is a complete-turn roll maneuver for the military aircraft Mirage-III, which is a single-seat, single-engine, fighter aircraft produced by the manufacturing company Dassault Aviation [21] for the French Air Force but widely exported and operated by Australia, Argentina, Pakistan, South Africa, Egypt [22] and other countries. The maximum take-off thrust is about 71 kN. Figure **9** shows three orthogonal views of the Mirage-III.

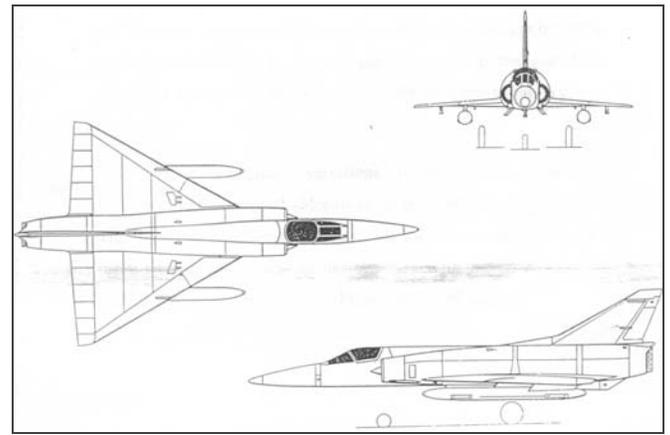

**Figure 9:** Drawings for Mirage-III.

The aircraft flies along a straight level path while performing a 360° continuous roll during a total time of 6 seconds with a constant velocity of 200 m/s at an altitude of 10 km. The speed of sound at that altitude is 299 m/s, so the Mach number is 0.67, which represents a high-speed subsonic flight without shock waves [23]. Thus, the aerodynamic treatment presented in the subsection (Aerodynamic and Stability Coefficients) is sufficiently adequate. In order to prescribe this maneuver, the constraint equations will be as follows:

$x_g(t) = 200\ t$

$y_g(t) = 0$

$z_g(t) = -10000$

$$\phi(t) = \left(\frac{2\pi}{16}\right)\left[\cos\left(\frac{3\pi t}{6}\right) - 9\cos\left(\frac{\pi t}{6}\right) + 8\right]$$

The profiles of the bank angle and its time derivatives are shown in Figure **10**. The maneuver duration is too short to worry about any change in the aircraft mass due to fuel consumption, and the constant-mass restriction in the mathematical model is of no concern. The air is assumed to be still (zero wind speed); therefore the relative aircraft/wind velocity reduces to the absolute velocity of the aircraft CG.

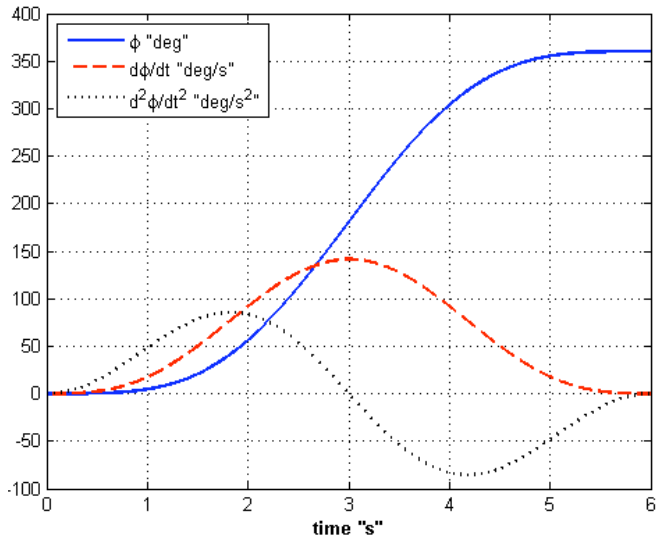

**Figure 10:** Temporal profile of the bank angle and its first and second derivatives.

### 8.1. Simulation Settings

The characteristic data used in the simulation are

#### 8.1.1. Mass, Inertia, and Main Dimensions

m = 7400 kg        A = 90000 kg.m$^2$        B = 54000 kg.m$^2$

C = 60000 kg.m$^2$    D = 0 kg.m$^2$           E = 1800 kg.m$^2$

F = 0 kg.m$^2$       d = b = 5.25 m          S = 36 m$^2$

#### 8.1.2. Aerodynamic and Longitudinal Stability Coefficients

$C_{L\alpha}$ = 2.204 rad$^{-1}$   $C_{L0}$ = 0.0[1]   $C_{D0}$ = 0.015   $K_{CD}$ = 0.4

$C_{m0}$ = 0.0    $C_{m\alpha}$ = –0.17 rad$^{-1}$    $C_{mq}$ = –0.4 rad$^{-1}$    $C_{m\delta m}$ = –0.45 rad$^{-1}$

#### 8.1.3. Lateral Stability Coefficients

$C_{y\beta}$ = –0.60 rad$^{-1}$

$C_{l\beta}$ = –0.05 rad$^{-1}$    $C_{lp}$ = –0.25 s/rad    $C_{lr}$ = 0.06 s/rad

$C_{l\delta l}$ = –0.30 rad$^{-1}$    $C_{l\delta n}$ = 0.018 rad$^{-1}$

$C_{n\beta}$ = 0.15 rad$^{-1}$    $C_{np}$ = 0.055 s/rad    $C_{nr}$ = –0.7 s/rad

$C_{n\delta l}$ = 0.0 rad$^{-1}$    $C_{n\delta n}$ = –0.085 rad$^{-1}$

#### 8.1.4. Flight Conditions

Altitude = 10000 m    $\rho^2$ = 0.412 kg/m$^3$    g = 9.81 m/s$^2$

With a time step of $\Delta t = 10^{-4}$ s (60001 time stations covering the 6-second maneuver), we obtained the following results from the computerized procedure. The solution was verified to be insensitive to the time step size because we used larger time steps of $2\times 10^{-4}$ s and $10^{-3}$ s and did not identify a notable difference. However, a too large time step like $10^{-2}$ s causes remarkable deviations with erratic profiles not capturing steep changes properly. The calculations take a few seconds on a laptop having a Core i5 CPU (4 cores) with 2.66 GHz speed. Referring to Equations (10-12), this maneuver has zero values for the flight path angles: $\theta_w$ and $\psi_w$. This was checked and found to be perfectly satisfied in our calculations.

### 8.2. Solution Results

The time-response required of the required thrust force is plotted in Figure **11**. The thrust profile is nearly symmetric in this maneuver. It experiences large and steep variations, but remains always positive and thus no reverse thrust is needed. During this maneuver, the angle of attack will drop and thus the lift coefficient will drop accordingly. The drag coefficient has a fixed part and another part that depends quadratically on the lift coefficient. With the lift coefficient dropping, the drag coefficient will drop sharply, resulting in sharp reduction in the demanded thrust for overcoming the drag force. The situation is inverted near the middle of the maneuver, where the angle of attack recovers its magnitude (but with a negative value), and thus the lift-induced drag steeply increases, elevating the demand on the thrust.

The necessary deflections angles of the rudder, elevator, and aileron during the maneuver are given in Figure **12**. The rudder shows the largest deflection, reaching a maximum of 49.9°. The expected variations of the actual angle of attack (as would be interpreted by an aerodynamicist, measured from the vertical $C_L$ axis in the $C_L$-$\alpha$ curve) and the sideslip angle are presented in Figure **13**. The angle of attack here is related to the one obtained from our solver by a constant shift of 6.36°. During the entire maneuver, the actual α lies between –6.05° and 6.36°. This should be a very safe range to apply our aerodynamic model that assumes a linear $C_L$-α relationship, being far below the stall value for the angle of attack, which is typically around 15°.

---

[1] According to subsection (About CL0 and $\alpha$) on page 19, this $C_{L0}$ is not actually needed here.
Instead, $C_{L0,equb}$ = (m g)/($\bar{q}$ S) = 0.245 is used. With $C_{L\alpha}$ = 2.204 rad$^{-1}$, we get $\alpha_{equb}$ = 0.111 rad (6.36°).

[2] Using the equation in subsection (Atmospheric Air Density).

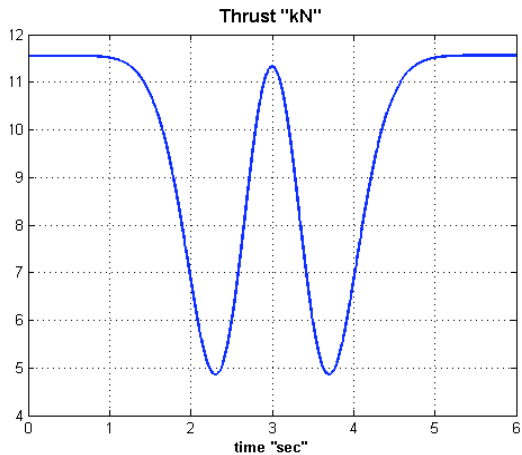

**Figure 11:** Time response of the thrust.

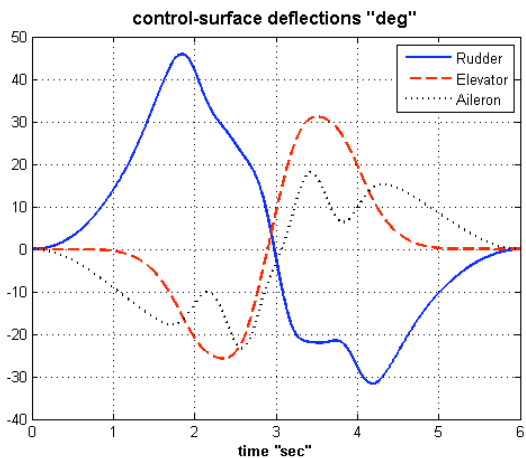

**Figure 12:** Time response of the deflection angles.

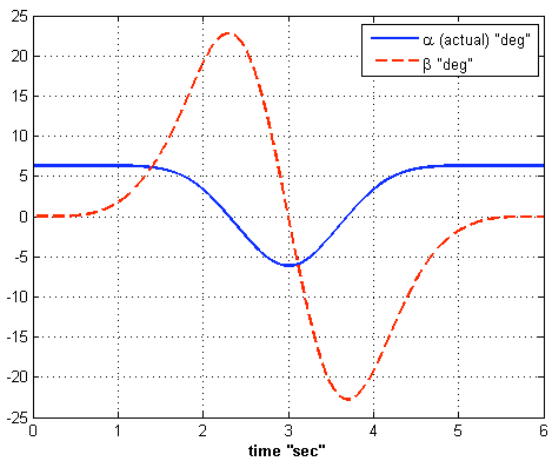

**Figure 13:** Time response of the actual angle of attack and the sideslip angle.

The responses of the pitch and yaw angles are presented in Figure **14**, and their phase portrait is depicted in Figure **15**. The former is nearly symmetric and always positive, and the latter is nearly anti-symmetric.

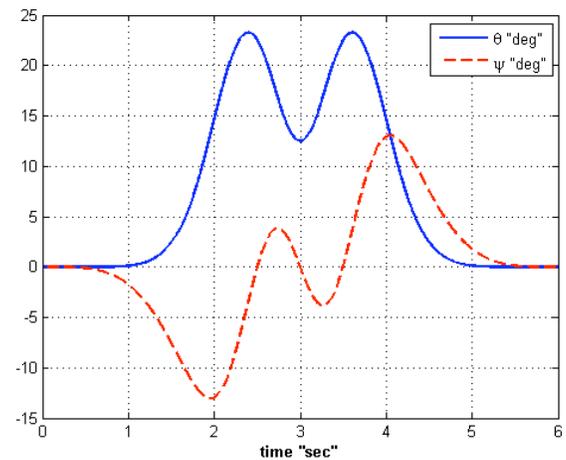

**Figure 14:** Time response of the pitch angle and the yaw (heading) angle.

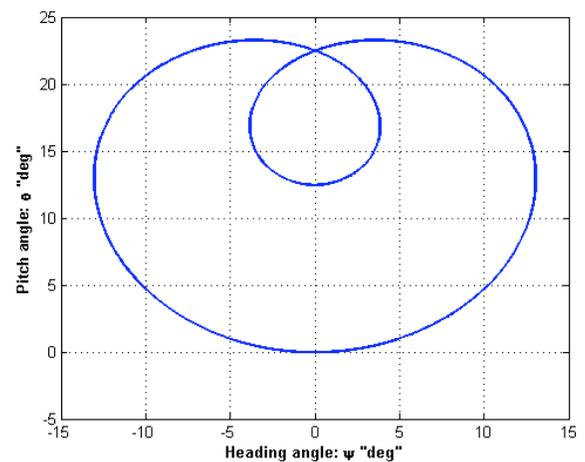

**Figure 15:** Phase portrait for the pitch angle versus the heading angle.

## 9. CONCLUSION

We presented a 3D nonlinear flight mechanics model for a generic fixed-wing aircraft, but having a constant mass and undergoing a shock-free stall-free flight. The model utilizes wind axes, Euler angles, local axes, and aerodynamic-performance expressions. The linear momentum equations were formulated in the wind axes to take advantage of the expected higher computational efficiency while the angular momentum equations were formulated in the body axes to take advantage of the time-invariance of the moments of inertia. Both sets of momentum equations are augmented with other kinematic, geometric, and aerodynamic relationships leading to a differential-algebraic system of 18 equations in 18 variables. Four out of these equations are actually provided as inputs describing the trajectory coordinates in the ground-fixed axes plus any other flight variable, such as the bank angle (as followed here). We presented in detail a computational algorithm which uses the 4$^{th}$ order Runge-Kutta integration method and time derivatives of the model equations to solve the nonlinear system sequentially in time to predict all the remaining 14 flight variables, especially the 3 deflection angles of the movable control surfaces (rudder, elevator, and ailerons) plus the magnitude of the thrust demanded from the engine. We implemented the algorithm using the software package MATLAB and applied it to a continuous roll maneuver of the fighter aircraft Mirage-III. The model and its implementation provide a useful design tool to quickly explore the upper limits of the thrust demand, control surface deflections, and angle of attack so that a proper aircraft design can be made taking these requirements into consideration. It can also be thought of as a pre-requisite for a programmed unmanned flight for a small aircraft where a microcontroller is installed in the aircraft and operates servo motors for deflecting the control surfaces and operates the engine throttle in such a way that achieves the predicted time responses of the control surface deflections and thrust for a limited period during a maneuver.

## ACKNOWLEDGEMENT

The author received cruicial support from Prof. Mohamed M. Abdelrahman (Department of Aerospace Engineering, Cairo University, Egypt), including the characterstic data of the Mirage-III and its maneuver case.

## REFERENCES

[1]   Avanzini G, Thomson D, Torasso A. Model predictive control scheme for rotorcraft inverse simulation. In *AIAA Guid Navig Control Conf*, Portland, Oregon, USA, 2011.

[2]   Zhou W, Wang H. Researches on inverse simulation's applications in teleoperation rendezvous and docking based on hyper-ellipsoidal restricted model predictive control for inverse simulation structure. J Aerosp Eng. 2014; 229(9): 1675-1689.

[3]   Lu L. Inverse modelling and inverse simulation for engineering applications, Saarbrücken, Germany: Lambert, 2010.

[4]   Blajer W, Graffstein J, Krawczyk M. Modeling of aircraft prescribed trajectory flight as an inverse simulation problem. In *Simulation and control of nonlinear engineering dynamical systems*, Netherlands, Springer, 2009; pp. 153-162.

http://dx.doi.org/10.1007/978-1-4020-8778-3_14

[5]   Leinonen A. Application for inverse simulation of flight tracks. *SCIENCE + TECHNOLOGY Aalto University publication series (Helsinki, Finland),* 2012; pp. 1-71.

[6]   Forster KJ, White TR. Numerical investigation into vortex generators on heavily cambered wings. AIAA J. 2014; 52(5): 1059-1071.

http://dx.doi.org/10.2514/1.J052529

[7]   Allerton D. Principles of flight simulation, United Kingdom: Wiley, 2009.

http://dx.doi.org/10.2514/4.867033

[8]   NASA Technical Reports Server (NTRS). Computer Mechanization of 6-DOF Flight Equations. NASA.

[9]   DAR Corporation: Design, Analysis and Research Corporation. AAA: Advanced Aircraft Analysis. [Online]. Available: www.darcorp.com/Software/AAA/. [Accessed 4 October 2015].

[10]  Anderson JD. Introduction to flight (5th ed), Indiana, USA: McGraw Hill Higher Education, 2005.

[11]  Abbott IH, Von Doenhoff AE. Theory of wing sections: including a summary of airfoil data, Dover Publications, 1959.

[12]  McClamroch NH. Steady aircraft flight and performance, USA: Princeton University Press, 2011.

[13]  Hoffmann KA, Chiang ST. Chapter 2: finite difference formulations. In *Computational fluid dynamics for engineers - volume I*, Wichita, Kansas, USA, Engineering Education System, 1993; pp. 30-60.

[14]  Mokhov II, Akperov MG. Tropospheric lapse rate and its relation to surface temperature from reanalysis data. Izvestiya, Atmos Oceanic Phys. 2006; 42(4): 430-438.

http://dx.doi.org/10.1134/S0001433806040037

[15]  Hibbeler RC. Engineering mechanics: Dynamics (13th ed), New Jersey: Pearson Prentice Hall, 2013.

[16]  Meriam JL, Kraige LG. Engineering mechanics - volume 2 dynamics, 7th ed, USA: Wiley, 2012.

[17]  Boeing. What is angle of attack? Aero Magazine 2000; 13.

[18]  Crawford B. Axes and derivatives. In *Flightlab Ground School*, Plymouth, MA, USA, Flight Emergency & Advanced Maneuvers Training, Inc, 2009; p. 1.2.

[19]  Fogarty LE, Howe RM. Exo Cruiser. 28 September 2011. [Online]. Available: http://dodlithr.blogspot.com/2011/09/computer-mechanization-of-6-dof-flight.html. [Accessed 31 October 2015].

[20]  MathWorks, "MATLAB and Simulink for Technical Computing," MathWorks, [Online]. Available: www.mathworks.com/. [Accessed 6 October 2015].

[21]  Dassault Aviation, "Dassault Aviation," Dassault Group, [Online]. Available: www.dassault-aviation.com/en/. [Accessed 6 October 2015].

[22]  Brindley J. Dassault mirage variants, Windsor, Berkshire, UK: Profile Publications, 1971.

[23]  Cole JD, Cook LP. Transonic aerodynamics, Amesterdam, The Netherlands: Elsevier Science Publishers, 2012.